\documentclass[10pt,twoside]{amsart}
\usepackage{amsmath}
\usepackage{amssymb}
\usepackage{amscd}
\usepackage{xypic}
\usepackage[all]{xy}

\title[Homology exponents for $H$-spaces]
{Homology exponents for $H$-spaces}

\author{Alain Cl\'ement}
\thanks{}

\address{
\hfill\break Rue Louis--Meyer 9\\
\hfill\break CH -- 1800 Vevey\\
\hfill\break Switzerland.} \email{alain.clement@bluewin.ch}

\author{J\'er\^ome Scherer}
\thanks{The second author is supported by the program Ram\'on y Cajal, MEC,
Spain, and FEDER/MEC grant MTM2004-06686. This research was
partially supported by the Swiss National Science Foundation grant
FN 200020-105383}
\address{
\hfill\break Departament de Matem\`atiques\\
\hfill\break Universitat Aut\`onoma de Barcelona\\
\hfill\break E--08193 Bellaterra\\
\hfill\break Spain.}
\email{jscherer@mat.uab.es}


\newcommand{\A}{{\mathcal A}}
\newcommand{\Ab}{\ifmmode{{\cal A}b}\else${\cal A}b$\fi}

\newcommand{\F}{{\mathbb F}}

\newcommand{\C}{{\mathbb C}}

\newcommand{\Z}{{\mathbb Z}}
\newcommand{\N}{{\mathbb N}}

\newcommand{\epi}{\twoheadrightarrow}

\newcommand{\degst}{\operatorname{deg}\nolimits}

\newtheorem{theorem}{Theorem}[section]
\newtheorem{proposition}[theorem]{Proposition}
\newtheorem{corollary}[theorem]{Corollary}
\newtheorem{lemma}[theorem]{Lemma}

\theoremstyle{definition}
\newtheorem{definition}[theorem]{Definition}
\newtheorem{remark}[theorem]{Remark}

\newtheorem{convention}[theorem]{Convention}
\newtheorem{notation}[theorem]{Notation}

\begin{document}

\begin{abstract}
We say that a space $X$ admits a \emph{homology exponent} if there
exists an exponent for the torsion subgroup of $H^*(X;\Z)$. Our
main result states if an $H$-space of finite type admits a
homology exponent, then either it is, up to $2$-completion, a
product of spaces of the form $B\Z/2^r$, $S^1$, $\C P^\infty$, and
$K(\Z,3)$, or it has infinitely many non-trivial homotopy groups
and $k$-invariants. We then show with the same methods that simply
connected $H$-spaces whose mod $2$ cohomology is finitely
generated as an algebra over the Steenrod algebra do not have
homology exponents, except products of mod $2$ finite $H$-spaces
with copies of $\C P^\infty$ and~$K(\Z,3)$.
\end{abstract}

\maketitle

\section*{Introduction}
The study of the torsion in the homotopy groups and the integral
homology groups of a space motivated the Moore conjecture, see
\cite{MR952582}, and the Serre conjecture, \cite{MR0060234}. Serre
proved that a simply connected space with finite dimensional (and
non-trivial) mod $p$ (co)homology $H^*(X; \F_p)$ must have
infinitely many non-trivial homotopy groups. He conjectured that
there should in fact exist infinitely many homotopy groups of $X$
containing $p$-torsion, which was proved eventually by McGibbon
and Neisendorfer \cite{MR749108}, relying on Miller's solution
\cite{Miller} of the Sullivan conjecture. This was then refined
further by Lannes and Schwartz in \cite{MR827370}. Their criterion
is that $H^*(X; \F_p)$ is locally finite, as a module over the
Steenrod algebra. Dwyer and Wilkerson went one step further,
\cite{MR92b:55004}, looking only at the module $QH^*(X; \F_p)$ of
indecomposable elements. F\'elix, Halperin, Lemaire, and Thomas
provided yet another criterion involving the depth of $H_*(\Omega
X; \F_p)$, \cite{MR974903}. In their subsequent paper
\cite{MR1184760} they focused on the size of the torsion part in
the ``loop space homology" $H_*(\Omega X; \Z)$. They proved in
fact a homological version of the Moore conjecture, namely that
the $p$-torsion part of the integral homology of the loop space of
a $\Z_{(p)}$-elliptic space always has an exponent.

In this article we are interested in understanding when the
torsion subgroup of the integral homology of a large class of loop
spaces, and more generally $H$-spaces, can have an exponent. In
the spirit of Serre's theorem, we first classify those $H$-spaces
having a homology exponent at the prime $2$ which are Postnikov
pieces (they have only a finite number of non-trivial homotopy
groups). Thus we will say henceforth that a space admits a
\emph{homology exponent} if there exists an integer $k$ such that
$2^k \cdot T_2 H^*(X; \Z) = 0$, where $T_2$ stands for the
$2$-torsion subgroup. We work with connected $H$-spaces of finite
type.



\medskip



\medskip

\noindent{\bf Theorem \ref{t:main}} {\it Let $X$ be an $H$-space
of finite type which admits a homology exponent. Then either $X$
is, up to $2$-completion, a product of spaces of the form
$B\Z/2^r$, $S^1$, $\C P^\infty$ and $K(\Z,3)$, or $X$ admits
infinitely many non-trivial $k$-invariants and homotopy groups.}

\medskip

The methods we develop predict in fact explicit degrees in which
to find homology classes of order $2^r$ for arbitrarily large $r$
when the space has no homology exponent, quite in the spirit of
Browder's ``infinite implications", \cite{MR23:A2201}. This builds
on previous work by the first author, who analyzed the case of a
Postnikov piece with at most two non-trivial homotopy groups in
\cite{MR2222506}.

There is a class of $H$-spaces which is very close to the
Postnikov pieces we have been dealing with up to now, namely those
$H$-spaces for which the mod $2$ cohomology is finitely generated
as an algebra over the Steenrod algebra. They are obtained indeed
as extensions by $H$-fibrations of an $H$-space with finite mod
$2$ cohomology by a Postnikov piece, \cite{deconstructing}.

\medskip

\noindent{\bf Theorem \ref{t:main2}} {\it Let $X$ be a simply
connected $H$-space of finite type such that $H^*(X; \F_2)$ is
finitely generated as an algebra over the Steenrod algebra. Assume
that $X$ admits a homology exponent. Then $X$ is, up to
$2$-completion, the product of a mod $2$ finite $H$-space $Y$ with
copies of $K(\Z, 2)$ and~$K(\Z,3)$.}

\medskip

This contrasts with the homological version of the Moore
conjecture obtained by F\'elix, Halperin, and Thomas in
\cite{MR1184760}. Of course the mod $2$ cohomology of the loop
space on a finite complex is very rarely finitely generated as an
algebra over the Steenrod algebra.

\medskip

\emph{Acknowledgements}. We would like to thank Richard Kane for
providing a simple proof of Lemma~\ref{lemma Qodd}, and Juan A.
Crespo and Wolfgang Pitsch for helpful comments. The second author
would like to thank Kathryn Hess and the IGAT, EPFL, for the
invitation which made this collaboration possible.

\section{Reduction to simply connected spaces}
\label{section reduction}

In this section we explain how to reduce the study of arbitrary
connected $H$-spaces to simply connected ones for which the second
homotopy group is torsion. These are then the spaces we study in
the rest of the article. Let us start with basic terminology and
notation.

\begin{notation}
\label{notation postnikov}
A space $X$ is a \emph{Postnikov piece} if it has only finitely
many non-trivial homotopy groups. It is an $H$-Postnikov piece if
it is moreover an $H$-space. The $n$-th \emph{Postnikov section}
$i_n: X \rightarrow X[n]$ is determined, up to homotopy, by the
property that it induces isomorphisms on homotopy groups $\pi_i$,
for $i \leq n$, and $\pi_i X[n] = 0$ for $i > n$. The homotopy
fiber $X \langle n \rangle$ of $i_n$ is the \emph{$n$-connected
cover} of $X$. When $X$ is simple (for example when $X$ is simply
connected or when $X$ is an $H$-space), there exist
\emph{$k$-invariants} $k_n \in H^{n+1}(X[n-1]; \pi_n X)$ such that
$X[n]$ can be recovered as the homotopy fiber of a map $k_n:
X[n-1] \rightarrow K(\pi_n X, n+1)$ representing the
$k$-invariant. When $X$ is an $H$-space, all $k$-invariants are
primitive elements.

Let $X$ be a space. By $\{B^*_r,d_r\}$ we denote its mod-$2$
cohomology \emph{Bockstein spectral sequence}: $B_1^*\cong
H^*(X;\F_2)\Longrightarrow (H^*(X;\Z)/\text{torsion})\otimes\F_2$.
Recall that the first differential $d_1 = Sq^1$ is the Bockstein
and a pair of elements $x$ and $y$ which survive to the page $B_r$
and such that $d_r(x) = y$ detect a copy of $\Z/2^r$ in $H^*(X;
\Z)$ in degree $|y| = |x| +1$.
\end{notation}

We collect now a result about ``small" Postnikov pieces. These
will turn out to be the only $H$-Postnikov pieces having an
exponent.

\begin{proposition}
\label{prop X[3]}
Let $X$ be a connected $H$-space of finite type such that
$\pi_2(X)$ and $\pi_3(X)$ are torsion free. Then $X[3]$ is a
product of spaces of the form $B\Z/p^r$, $S^1$, $\C P^\infty$ and
$K(\Z,3)$. Moreover, $X[3]$ admits a homology exponent.
\end{proposition}

\begin{proof}
It is well-known that a copy of the integers in $\pi_1 X$
corresponds to a copy of $S^1$ splitting of $X$ (because of the
existence of a section $S^1 \rightarrow X$). One readily verifies
that $H^3(K(\Z/p^r,1);\Z)=0$, which shows that the first
$k$-invariant must be trivial. Thus $X[2]$ splits as a product of
copies of $B\Z/p^r$'s, $S^1$'s, and $\C P^\infty$'s. Next, the
only elements in $H^4(K(\Z/p^r,1);\Z)$, for any prime $p$ and any
integer $r$, and $H^4(K(\Z,1);\Z)$ are multiples of the square of
the generator in degree~$2$. Such elements are not primitive
(unless they are trivial), and hence cannot be the $k$-invariants
of an $H$-space. Therefore, the second $k$-invariant of $X[3]$ is
trivial as well and the space splits as a product.

It remains to prove the assertion about the homology exponent.
Recall that $p^r\cdot\widetilde{H}^*(B\Z/p^r;\Z)=0$ (by a transfer
argument), $H^*(S^1;\Z)$, $H^*(\C P^\infty;\Z)$ are torsion free
and $2\cdot T_2(H^*(K(\Z,3);\Z))=0$ (as a consequence of Serre
\cite{MR0060234} or Cartan's computations \cite{MR0087934}).
\end{proof}

\begin{lemma}
\label{lemma circlefiber}
Let $r\geq1$ and $(S^1)^r\to Y\to X$ be a $H$-fibration. If $X$
admits a homology exponent, then so will $Y$.
\end{lemma}

\begin{proof}
Let us remark that the fibration is orientable \cite[p.
476]{MR666554}. We obtain the result for $r=1$ by inspecting the
associated Gysin cohomology exact sequence and conclude by
induction on~$r$.
\end{proof}

We conclude this section with the promised reduction. The
existence of a homology exponent for an arbitrary $H$-space is
detected in the homology of a certain covering space.

\begin{proposition}
\label{prop reduction}
Let $X$ be a connected H-space of finite type. Then there exists a
simply connected $H$-space of finite type $Y$ such that $\pi_2Y$
is a torsion abelian group and $Y$ fits into the following
$H$-fibration:
$$\xymatrix{
Y\ar[r] &X\ar[r] &B\pi_1X\times K(\Z^k,2), }$$
for some $k\geq0$.
Moreover, if $X$ admits a homology exponent, then so will $Y$.
\end{proposition}

\begin{proof}
Let us first deal with the copies of $\Z$ in $\pi_1 X$. They
correspond to a torus $(S^1)^r$ splitting off $X$. There exists
hence $X_1$ such that $X\simeq X_1\times(S^1)^r$ and $\pi_1
X_1\cong\oplus_a\Z/2^{s_a} \oplus A$, where $A$ is a $2'$-torsion
abelian group. For any $a$ there is a map $B\Z/2^{s_a}\to K(\Z,2)$
corresponding to the Bockstein operation of order $s_a$. Let us
define $X_2$ to be the homotopy fiber of the composite map
$$
X_1\to B\pi_1X_1\to \prod_a K(\Z,2) \times K(A, 1).
$$
As in the proof of \cite[Proposition 0.7]{MR827370}, $X_2$ splits
as a product $\widetilde{X}\times\prod_a S^1$, where
$\widetilde{X}$ denotes the universal cover of $X$. By the
previous lemma, $X_2$, and thus $\widetilde{X}$, admit an exponent
if $X$ does.

Finally let us write $\pi_2 X \cong \Z^k \oplus A'$ where $A'$ is
a finite torsion group and define $Y$ to be the homotopy fiber of
the map $\widetilde{X}\rightarrow K(\Z^k, 2)$ so that $\pi_2 Y
\cong A'$. The previous lemma yields the statement about the
exponent. The description of the base of the $H$-fibration comes
from Proposition~\ref{prop X[3]}.
\end{proof}

\section{A splitting principle}
\label{section splitting}

Let $X$ be a Postnikov piece, which highest non-trivial homotopy
group is $\pi_n X = H$. In this section we show that, if $H$ is
torsion free, the $(n-2)$-connected cover $X(n-2)$ splits as a
product $K(H, n) \times K(G, n-1)$. Loosely speaking, the $n$-th
$k$-invariant attaches $K(H, n)$ directly to $X[n-2]$ and cannot
tie the last two homotopy groups together.

We will first need some basic results on Eilenberg-Mac\,Lane
spaces. We follow the terminology and notation of
\cite[Chapter~1]{MR95d:55017}. For a unified treatment of the
spaces $K(\Z/2^s, n)$, with $s \geq 1$, and $K(\Z, n)$, it is
convenient to introduce a notation for the higher Bockstein
operations. Let $u_n$ (respectively $Sq^1_s u_n$) denote the
generator of the $1$-dimensional $\F_2$-vector space
$H^n(K(\Z/2^s,n);\F_2)$ (respectively
$H^{n+1}(K(\Z/2^s,n);\F_2)$). For an admissible sequence $I =
(i_1, \dots, i_m)$, we will write $Sq^I_s u_n$ instead of
$Sq^{(i_1, \dots, i_{m-1})}Sq^1_s u_n$ if $i_m =1$ and instead of
$Sq^I u_n$ if $i_m \neq 1$. We denote also by $u_n$ the generator
of $H^n(K(\Z,n);\F_2)$.

Serre computed the mod-$2$ cohomology of Eilenberg-Mac\,Lane
spaces.

\begin{theorem} {\rm (Serre, \cite{MR0060234})}
\label{theorem Serre}
Let $n\geq1$ and $s\geq1$. \begin{itemize}
\item[(1)] The $\F_2$-algebra
$H^*(K(\Z/2^s,n);\F_2)$ is isomorphic to the polynomial
$\F_2$-algebra on generators $Sq^I_s u_n$, where $I$ covers all
the admissible sequences of excess $e(I)<n$.
\item[(2)] The $\F_2$-algebra
$H^*(K(\Z,n);\F_2)$ is isomorphic to the polynomial $\F_2$-algebra
on generators $Sq^I u_n$, where $I$ covers all the admissible
sequences of the form $(i_1,\dots,i_n)$ where $i_n \neq 1$ and of
excess $e(I)<n$. \hfill{\qed}
\end{itemize}
\end{theorem}

Our first lemma relies on Serre's computations in low degrees.

\begin{lemma}
\label{lemme [K(G,n-1),K(H,n+1)]}
Let $G$ be a finitely generated abelian group, $H$ be free
abelian, and $n\geq3$. Then $P^{n+1}H^*(K(G,n-1);H)=0$.
\end{lemma}

\begin{proof}
Since $H^*(K(G,n-1);H) \cong H^*(K(G,n-1); \Z) \otimes H$, it is
enough to consider the case when $H = \Z$. When $n \geq 4$, the
only elements in $H^{n+1}(K(G,n-1); \F_2)$ are sums of elements of
the form $Sq^2 u_{n-1}$. These elements all have non-trivial
Bockstein and we see therefore from the Bockstein spectral
sequence that $H^{n+1}(K(G,n-1); \Z)$ is $2$-torsion free. As
there is obviously no odd primary torsion in this degree (from
Cartan's description of $H^*(K(G, n-1); \F_p)$, \cite{MR0087934}),
we have that $H^{n+1}(K(G,n-1); \Z)=0$. When $n =3$, write $G
\cong \Z^k \oplus A$, where $A$ is torsion. Then $H^{4}(K(G,2);
\Z) \cong H^{4}(K(\Z^k,2); \Z)$. There are no primitive elements
in this degree.
\end{proof}

\begin{proposition}
\label{proposition factorisation du k-invariant}
Let $X$ be a simply connected $H$-space of finite type such that
$\pi_2X$ is a torsion group. Let $n\geq3$, consider the Postnikov
section $X[n]$ and assume that $\pi_nX$ is torsion free. Then we
have the following $H$-fibration:
$$\xymatrix{
X[n]\ar[r] &X[n-2]\ar[r] &K(\pi_{n-1}X,n)\times K(\pi_nX,n+1).
}$$
\end{proposition}

\begin{proof}
Set $G=\pi_{n-1} X$. We prove that the $k$-invariant $k_n \in
H^{n+1}(X[n-1]; H)$ factors through $H^{n+1}(X[n-2]; H)$. Let us
consider the fibration
$$
K(G, n-1) \rightarrow X[n-1] \rightarrow
X[n-2]
$$
and the cofibration $K(G, n-1) \rightarrow X[n-1]
\rightarrow C$. From Lemma~\ref{lemme [K(G,n-1),K(H,n+1)]} we
deduce that the $H$-map $K(G, n-1) \rightarrow X[n-1]
\xrightarrow{k_n} K(H, n+1)$ is null-homotopic.

Therefore $k_n$ factors through a map $C \rightarrow K(H, n+1)$.
By Ganea's result \cite{MR31:4033} the fiber of the map $C
\rightarrow X[n-2]$ is the join $K(G, n-1)*\Omega (X[n-2])$. This
is an $n$-connected space and $\pi_{n+1} (K(G, n-1)*\Omega
(X[n-2])) \cong G \otimes \pi_2 X$, a torsion group because $\pi_2
X$ is so. Thus $K(H, n+1)$ is a $K(G, n-1)*\Omega (X[n-2])$-local
space, as we assume that $H$ is torsion free. From Dwyer's version
of Zabrodsky lemma \cite[Proposition~3.5]{MR97i:55028} we deduce
that $k_n$ factors through a map $k: X[n-2] \rightarrow K(H,
n+1)$.

If $k_{n-1} \in H^n(X[n-2]; G)$ denotes the previous
$k$-invariant, this means that $X[n]$ is the homotopy fiber of the
product map $X[n-2] \xrightarrow{k_{n-1} \times k} K(G,n)\times
K(H,n+1)$.
\end{proof}

\section{Gaps in the primitives}
\label{section gaps}

This section contains the key cohomological result which makes the
analysis of the Serre spectral sequence possible. We notice first
that there are gaps in the mod $2$ cohomology of
Eilenberg-Mac\,Lane spaces and show then that these gaps propagate
in the cohomology of any Postnikov piece.

\begin{definition}
\label{A_n}
Let $n\geq1$. We set $A_n=\{a\in\N\text{, $a$ odd}\ |\
\nu_2(a)\geq n+1\}$ where $\nu_2(a)$ is the $2$-adic length of the
integer $a$.
\end{definition}

We will show that there are no indecomposable elements in the
cohomology of an $n$-stage Postnikov piece in degrees $a \in A_n$.
To deduce that there are no primitive elements either, we make use
of the relationship provided by the Milnor-Moore theorem
\cite[Proposition 4.21]{MR0174052}: For a connected, associative,
and commutative Hopf algebra over $\F_2$, there is an exact
sequence of graded modules
$$
\xymatrix{ 0\ar[r] &P(\xi H)\ar[r]
&PH\ar[r] &QH, }
$$
where $\xi H$ is the image of the \emph{Frobenius map}
$\xi:x\mapsto x^2$, $QH$ is the module of indecomposable elements
and $PH$ is the module of primitive elements of~$H$.

\begin{lemma}
\label{lemme QK(H,n)}
Let $H$ be a finitely generated abelian group and $n\geq2$. Then
$$
Q^a H^*(K(H,n);\F_2)=0=P^a H^*(K(H,n);\F_2)
$$
for all $a\in A_n$.
\end{lemma}

\begin{proof}
When it is not trivial, the $\F_2$-algebra structure of
$H^*(K(H,n);\F_2)$ is given by a polynomial algebra on generators
of the form $Sq_s^Iu_n$ where $I$ runs over admissible sequences
with excess $e(I)<n$, as we have recalled in Theorem~\ref{theorem
Serre}. Careful calculations show that these generators lie in
degrees $1+2^{h_1}+\dots+2^{h_{n-1}}$ where $h_1\geq\dots\geq
h_{n-1}\geq0$ (see \cite[Th\'eor\`eme 1, p. 212 and Th\'eor\`eme
2, p. 213]{MR0060234}). The $2$-adic length of these degrees is
bounded by~$n$. This shows that there are no indecomposable
elements in the degrees we claimed. These degrees being odd, there
are no primitives either, because the kernel of the map
$PH^*(K(H,n);\F_2) \rightarrow QH^*(K(H,n);\F_2)$ is concentrated
in even degrees.
\end{proof}

The proof of the following lemma has been kindly communicated to
us by Richard Kane, \cite{Kanepers}.

\begin{lemma}
\label{lemma Qodd}
Let $B$ be a connected, associative, and commutative Hopf algebra
of finite type over $\F_2$ and $A$ a sub-Hopf algebra of~$B$. Then
the morphism $QA \rightarrow QB$ is injective in odd degrees.
\end{lemma}

\begin{proof}
We work in degree $2n+1$. Consider the Hopf subalgebra $C$ of $A$,
and hence of $B$, generated by the elements in $A$ of degree $\leq
2n$. Then one has an inclusion of quotient Hopf algebras $A//C
\hookrightarrow B//C$ by \cite[Corollary p.9]{Kane}. Let $x$ be an
indecomposable element in $QA$ of degree $2n+1$. It determines a
non-zero primitive element in $P(A//C)$, hence in $P(B//C)$. As
the map $P(B//C) \rightarrow Q(B//C)$ is injective in odd degrees,
we see that the composite $QA \rightarrow Q(B//C)$ is injective in
degree $2n+1$. Therefore $QA \rightarrow QB$ must be injective in
degree $2n+1$ as well.
\end{proof}

\begin{remark}
\label{remark andrequillen}
The preceding lemma has a nice interpretation in terms of
Andr\'e-Quillen homology, the derived functor of $Q(-)$. It is
proved in \cite[Proposition~1.3]{CCS5} that one has, in the
setting of the lemma, an exact sequence $H^Q_1(B//A)\rightarrow QA
\rightarrow QB \rightarrow Q(B//A)\rightarrow 0$, a result dual to
that of Bousfield, \cite[Theorem 3.6]{MR0258919}. Moreover the
graded $\F_2$-vector space $H^Q_1(B//A)$ is concentrated in even
degrees.
\end{remark}

We are now ready to prove that the gaps also appear in the
cohomology of any Postnikov piece.

\begin{proposition}
\label{proposition QX[n]}
Let $n\geq2$ and $X$ be a simply connected $n$-stage $H$-Postnikov
piece of finite type. Then $Q^a H^*(X;\F_2)= 0 = P^a H^*(X;\F_2)$
for all $a\in A_n$.
\end{proposition}

\begin{proof}
The proof goes by induction on $n$. We have the following
$H$-fibration given by the Postnikov tower of $X$:
$$
\xymatrix{ K(\pi_n(X),n)\ar[r]&X\ar[r]^-p&X[n-1]. }
$$
We rely on the analysis of the Eilenberg-Moore spectral sequence
done by Smith \cite[Proposition~3.2]{MR0275435}. The quotient Hopf
algebra $R=H^*(X[n-1];\F_2)//\ker p^*$ can be identified via $p^*$
with a sub-Hopf algebra of $H^*(X;\F_2)$. The corresponding
quotient $S=H^*(X;\F_2)//R$ is isomorphic to a sub-Hopf algebra
(and a sub-$\A_2$-algebra) of $H=H^*(K(\pi_n(X),n);\F_2)$. There
is a section $S \rightarrow H^*(X;\F_2)$, which is a map of
algebras, so that the module of indecomposables $QH^*(X;\F_2)$ is
isomorphic to $QR \oplus QS$, as graded $\F_2$-vector spaces. We
have to prove that both $QR$ and $QS$ are trivial in degrees in
$A_n$. First, since $Q(-)$ is right exact, we have a surjection
$QH^*(X[n-1];\F_2)\epi QR$. Now, $Q^aH^*(X[n-1];\F_2)=0$ for any
$a \in A_{n-1}$ by induction hypothesis and we conclude that
$Q^aR=0$ for any $a \in A_n$ since $A_n\subset A_{n-1}$. Second,
we deal with $QS$. Let us apply the preceding lemma to the
inclusion $S \subset H^*(X;\F_2)$. We see that $QS \rightarrow QH$
is a monomorphism in odd degrees. Therefore $Q^aS=0$ for all $a\in
A_n$ by Lemma~\ref{lemme QK(H,n)}.
\end{proof}

\section{Transverse elements in Eilenberg-Mac\,Lane spaces}
\label{section EML}

Now begins the study of the $2$-torsion in Postnikov pieces. In
this section we deal with the first step of the induction, namely
the analysis of the case of Eilenberg-Mac\,Lane spaces. Recall
that $\{B^*_r,d_r\}$ denote the mod-$2$ cohomology Bockstein
spectral sequence of a space~$X$.

\begin{definition}
\label{def transverse}
Let $n$ and $r$ be two positive integers. An element $x\in B^n_r$
is said to be \emph{$\ell$-transverse} if $d_{r+l}x^{2^l}\not=0\in
B^{2^l n}_{r+l}$ for all $0\leq l\leq\ell$. An element $x\in
B^n_r$ is said to be \emph{transverse} if it is $\ell$-transverse
for all $\ell\geq0$. We will also speak of \emph{transverse
implications} of an element $x\in B^n_r$.
\end{definition}

Every transverse element gives rise to $2$-torsion of arbitrarily
high order in the integral cohomology of $X$. This definition,
introduced in \cite{Cl02-PhD}, adapts Browder's ``infinite
implications" from \cite{MR23:A2201} to our purpose. To us, the
fact that the elements die in increasing pages of the Bockstein
spectral sequence is crucial, whereas Browder was merely
interested to know that the degrees of the elements was
increasing.

Our strategy for disproving the existence of a homology exponent
for a space will consist in exhibiting a transverse element in its
mod-$2$ cohomology Bockstein spectral sequence. Note that in
principle the absence of transverse elements does not imply the
existence of an exponent for the $2$-torsion part in $H^*(X;\Z)$.
An easy example if given by the infinite wedge $\vee M(\Z/2^n,
n)$.

In the special case of Eilenberg-Mac\,Lane spaces, we have the
following result, taken from the first author PhD thesis,
\cite[Theorem 1.3.2]{Cl02-PhD}.

\begin{theorem}
\label{t:transversity}
Let $H$ be an abelian group of finite type and let $n\geq2$.
Consider the Eilenberg-Mac\,Lane space $K(H,n)$ and its mod-$2$
cohomology Bockstein spectral sequence $\{B^*_r,d_r\}$. Suppose
that one of the following assumptions holds:

\begin{itemize}
\item[$\bullet$]{$n$ is even and $x\in B^n_{s_j}$ is $0$-transverse
for some $1\leq j\leq l$,}
\item[$\bullet$]{$x\in P^{\text{even}} B_1^*$ is $0$-transverse
($Sq^1x\not=0$).}
\end{itemize}
Then $x$ is transverse. \hfill{\qed}
\end{theorem}

Note that the abelian group $H$ is isomorphic to
$\Z^{s}\oplus\Z/2^{s_1}\oplus\dots\oplus\Z/2^{s_l} \oplus A$,
where $A$ is a $2'$-torsion group, which is therefore invisible to
the mod $2$ Bockstein spectral sequence. Hence the first type of
$0$-transverse elements correspond basically to the fundamental
classes $u_n$ introduced in Section~\ref{section splitting}, one
for each copy of $\Z/2^{s_j}$ (the fundamental classes coming from
the copies of $\Z$ survive to $B^n_\infty$).

\begin{remark}
\label{rem:transversity}
In general a $0$-transverse implication does not imply transverse
implications. More precisely, the fact that $x\in P^{\text{even}}
H^*(X;\F_2)$ is such that $Sq^1x\not=0$ does not always force $x$
to be transverse. A counter-example is given by $X=BSO$ and
$x=w_2$, the second Stiefel-Withney class in $H^2(BSO;\F_2)$.
\end{remark}

From Theorem \ref{t:transversity} it is not difficult to prove
that most Eilenberg-Mac\,Lane spaces have no homology exponent.

\begin{proposition}
\label{c:GEM_no_exp}
Let $H$ be a non-trivial $2$-torsion abelian group and let
$n\geq2$. The Eilenberg-Mac\,Lane space $K(H,n)$ has no homology
exponent.
\end{proposition}

\begin{proof}
Accordingly to the K\"unneth formula, it is sufficient to
establish the result when $H=\Z/2^s$ for some $s\geq1$. If $n$ is
even, consider the reduction of the fundamental class $u_n\in
H^n(K(\Z/2^s,n);\F_2)$. This class survives to $B^n_s$ and is
$0$-transverse. Then $u_n\in B^n_s$ is transverse. If $n$ is odd,
consider the admissible sequence $(2,1)$. Its excess is exactly
$1$ and therefore $Sq^{2,1}_s u_n\in P^{\text{even}}
H^*(K(\Z/2^s,n);\F_2)$ when $n\geq3$. Moreover we have
$Sq^1Sq^{2,1}_s u_n=Sq^{3,1}_s u_n$ by Adem relations, which means
that $Sq^{2,1}_s u_n$ is $0$-transverse. Hence $Sq^{2,1}_s u_n\in
B^{n+3}_1$ is transverse.
\end{proof}

\begin{proposition}
\label{prop KZn}
Let $H$ be a finitely generated abelian group and $n\geq4$. The
Eilenberg-Mac\,Lane space $K(H,n)$ is then either mod $2$ acyclic,
or has no homology exponent.
\end{proposition}

\begin{proof}
By the K\"unneth formula and Proposition \ref{c:GEM_no_exp}, it is
sufficient to analyze the case $H=\Z$. Consider the reduction of
the fundamental class $u_n\in H^n(K(\Z,n);\F_2)$. If $n$ is even,
then $Sq^{2}u_n$ is transverse. If $n$ is odd, then $Sq^{6,3}u_n$
is transverse.
\end{proof}

\section{Transverse elements in Postnikov pieces}
\label{section postnikov pieces}

We are now ready to prove our main result: Most Postnikov pieces
do not have a homology exponent. The strategy to prove this relies
on the crucial observation that the transverse implications of
certain element in the cohomology of the total space of a
fibration can be read in the cohomology of the fibre.

\begin{lemma}
Let $j:F\to X$ be a continuous map. If $x\in H^*(X;\F_2)$ is such
that $j^*(x)\not=0\in H^*(F;\F_2)$ is transverse, then $x$ itself
is transverse.
\end{lemma}

\begin{proof}
It follows from the naturality of the Bockstein spectral sequence.
\end{proof}

\begin{theorem}
\label{t:main}
Let $X$ be an $H$-space of finite type which admits a homology
exponent. Then either $X$ is, up to $2$-completion, a product of
spaces of the form $B\Z/2^r$, $S^1$, $\C P^\infty$ and $K(\Z,3)$,
or $X$ admits infinitely many non-trivial $k$-invariants and
homotopy groups.
\end{theorem}

\begin{proof}
Let us assume that $X$ is a Postnikov piece. By Proposition
\ref{prop reduction}, there is an $H$-fibration of the form
$$\xymatrix{
Y\ar[r] &X\ar[r] &B\pi_1X\times K(\Z^r,2), }$$ where $Y$ is a
simply connected $H$-space of finite type such that $\pi_2Y$ is a
torsion abelian group. Moreover, $Y$ admits a homology exponent.
It is also clearly a Postnikov piece. Let us show that $Y$ is, up
to $2$-completion, a product of copies of $K(\Z,3)$. By
Proposition \ref{prop X[3]}, this will imply that $X$ itself
splits as the announced product.

Assume that $\pi_nY = H$ is the highest non-trivial homotopy group
of~$Y$, up to $2$-completion. If $n=2$, since $\pi_2 Y$ is a
torsion abelian group, we deduce from
Proposition~\ref{c:GEM_no_exp} that $H$ is $2'$-torsion. In other
words $Y^\wedge_2$ is contractible. We can therefore assume that
$n \geq 3$. The space $Y$ fits into the fibration sequence
$$\xymatrix{
K(H,n)\ar[r]^-j &Y\ar[r]^-i &Y[n-1]\ar[r]^-k &K(H,n+1), }
$$
where $k$ denotes the last $k$-invariant. We analyze the situation
in two steps, depending on the presence of $2$-torsion in $H$.

Let us first assume that $H$ contains $2$-torsion, let us say
$\oplus_b\Z/2^{t_b}$. Choose an index $b$ and consider the
projection $\pi: H \rightarrow \Z/2^{t_b}$ on the corresponding
cyclic subgroup. Pick $v_n \in H^n(K(H, n); \F_2)$, the image via
$\pi^*$ of the class $u_n \in H^n(K(\Z/2^{t_b}, n); \F_2)$.

Set $\xi=(2^{n-1}-2,2^{n-2}-1,2^{n-3}-1,\dots,3,1)$. The degree
$\deg(Sq_t^\xi v_n)=2^n-2$ is even and $Sq^1Sq_t^\xi v_n\not=0$
since $e(\xi)=n-2$. By Theorem \ref{t:transversity}, $Sq_t^\xi
v_n$ is transverse. Since $Y$ is an $H$-space and the
$k$-invariant is an $H$-map, the element $d_{n+1}v_n$ is
primitive, and so is $d_{2^n-1}Sq_t^\xi v_n=Sq_t^\xi d_{n+1}v_n\in
P^{2^n-1}H^*(Y[n-1];\F_2)$. By Proposition \ref{proposition
QX[n]}, $P^{2^n-1}H^*(Y[n-1];\F_2)=0$ since $2^n-1\in A_{n-1}$.
Therefore, $Sq_t^\xi v_n$ survives in the Serre spectral sequence
and by the previous lemma, $H^*(Y;\F_2)$ contains a transverse
element. In particular it has no homology exponent.

Hence, $H$ must be $2$-torsion free and is thus isomorphic to
$\Z^s \oplus A$, where $A$ is a torsion group, for some $s \geq
1$. By Proposition \ref{proposition factorisation du k-invariant},
$Y$ fits in the following $H$-fibration:
$$\xymatrix{
K(H,n)\times K(\pi_{n-1}Y,n-1)\ar[r] &Y\ar[r] &Y[n-2].
}
$$
Choose now $v_n \in H^n(K(H, n); \F_2)$ to be the image of the
class $u_n \in H^n(K(\Z, n); \F_2)$ given by projection on the
first copy of $\Z$ in $H$.

If $n\geq4$, then set
$\eta=(2^{n-2}+2^{n-3}-2,2^{n-3}+2^{n-4}-1,2^{n-4}+2^{n-5}-1,\dots,5,2)$.
The degree $\deg(Sq^\eta v_n)=2^{n-1}+2^{n-2}-2$ is even and
$Sq^1Sq^\eta v_n\not=0$ since $e(\eta)=n-2$. Thus $Sq^\eta v_n$ is
transverse and survives in the Serre spectral sequence of the
above fibration since $2^{n-1}+2^{n-2}-1\in A_{n-2}$. In this
case, $H^*(Y;\F_2)$ contains a transverse element and has no
homology exponent.

Therefore, $n=3$ and $Y\simeq K(H,3)\times K(\pi_2Y,2)$. Since $Y$
admits a homology exponent, the torsion group $\pi_2Y$ is trivial
and $H$ is torsion free.
\end{proof}

The proof of the theorem predicts explicit degrees in which to
find higher and higher torsion in the integral cohomology of the
space.

\begin{corollary}
\label{explicit degrees}
Let $X$ be a simply connected $H$-Postnikov piece of finite type,
say $X \simeq X[n]$. Assume that $\pi_2 X$ is torsion and that $X$
is not equivalent up to $2$-completion to a product of copies of
$K(\Z, 3)$. Then, for any integer $k$, there is a copy of $\Z/2^k$
in $H^{*}(X; \Z)$
\begin{itemize}
\item[(1)] in degree $2^k(2^n-2)$ if $\pi_n X$ contains
$2$-torsion,
\item[(2)]in degree $2^k(2^{n-1}+2{n-2} -2)$ if not.
\end{itemize}
\end{corollary}

\begin{proof}
Since $X$ is different from $K(\Z^m, 3)$, we know from
Theorem~\ref{t:main} that $X$ has no exponent. The higher and
higher torsion is detected by the consecutive powers of the
elements $Sq^\xi v_n$ and $Sq^\eta v_n$ constructed in the above
proof.
\end{proof}

Any finite $H$-space has obviously a homology exponent. Our second
corollary applies to its Postnikov sections. As soon as it has at
least two homotopy groups, it cannot have a homology exponent.

\begin{corollary}
Let $X$ be a simply connected finite $H$-space and $n\geq3$. Then
$X[n]$ has a homology exponent if and only if $X[n]\simeq X[3]
\simeq K(\Z^r,3)$ for some $r\geq0$.
\end{corollary}

\begin{proof}
The fact that the $H$-space $X$ is finite and simply connected
forces it to be $2$-connected, \cite[Theorem~6.10]{MR23:A2201}.
Moreover, $\pi_3 X \cong \Z^r$ for some integer~$r$, by work of
Hubbuck and Kane, \cite{MR53:1582}. The result now follows
directly from Theorem~\ref{t:main}.
\end{proof}

This corollary applies in particular to $S^3$. The Postnikov
section $S^3[3] \simeq K(\Z, 3)$ has a homology exponent, but all
higher Postnikov sections $X[n]$, $n\geq4$, have none. The
following proof of a result obtained by Levi in \cite{MR1308466}
is, to our knowledge, the first one not based on Miller's solution
of the Sullivan's conjecture \cite{Miller}. Let us mention in this
context the work of Klaus, \cite{MR1879944}, who proves the
statement about the $k$-invariants for $BG^{\wedge}_2$, not for
the loop space.

\begin{corollary}
Let $G$ be a $2$-perfect finite group. Then $\Omega
(BG^{\wedge}_2)$ has infinitely many non-trivial $k$-invariants
and homotopy groups.
\end{corollary}

\begin{proof}
Suppose $BG^{\wedge}_2$ is a Postnikov piece. Following Levi
\cite{MR1308466}, there is a homology exponent for $\Omega
(BG^\wedge_2)$ and therefore this space has to be a product of
copies of $B\Z/2^r$, $S^1$, $\C P^\infty$ or $K(\Z,3)$. Since
$\Omega (BG^\wedge_2)$ has torsion homotopy groups, the only
copies that can occur are of the form $B\Z/2^r$. Thus
$BG^\wedge_2\simeq K(A,2)$, where $A$ is a $2$-torsion abelian
group. By the Evens-Venkov theorem, \cite{MR0137742},
$H^*(BG^\wedge_2;\F_2)$ is Noetherian. Hence $A$ is trivial, and
so is $BG^\wedge_2$.
\end{proof}

\section{Comparison with other forms of Serre's theorem}
\label{section compare}
In this section we compare our theorem to the other results we
mentioned in the introduction. We show that the existence of a
homology exponent is stronger than all previously established
criteria, except possibly \cite{MR974903}, which seems difficult
to relate directly to cohomological statements. Therefore, when
$X$ is an $H$-space, our result provides new proofs of those. They
are very different in spirit, since they do not require the
Sullivan conjecture. For simplicity we deal here with simply
connected spaces.

\begin{proposition}
Let $X$ be a simply connected $H$-Postnikov piece. Then
\begin{itemize}
\item[(1)] {\rm (Serre \cite{MR0060234})} $H^*(X; \F_2)$ is not finite,
\item[(2)] {\rm (Lannes-Schwartz \cite{MR827370})} $H^*(X; \F_2)$ is not locally finite,
\item[(3)] there exists an element of infinite height in $H^*(X;
\F_2)$,
\item[(4)] {\rm (Grodal \cite{MR1622342})} the transcendence degree of
$H^*(X; \F_2)$ is infinite unless $X$ is homotopy equivalent, up
to $2$-completion, to $K(\Z, 2)^s$,
\item[(5)] {\rm (Dwyer-Wilkerson \cite{MR92b:55004})} the unstable
module $QH^*(X; \F_2)$ is not locally finite unless $X$ is
homotopy equivalent, up to $2$-completion, to $K(\Z, 2)^s$.
\end{itemize}
\end{proposition}

\begin{proof}
Notice first that $K(\Z, 2)$ and $K(\Z, 3)$ satisfy (1) - (5).
Assume now that $X$ is a Postnikov piece, say $X \simeq X[n]$. In
the proof of Theorem~\ref{t:main} we first considered the covering
fibration $(S^1)^r \rightarrow Y \rightarrow X$. The map $Y
\rightarrow X$ induces isomorphisms in homology in high degrees.
We can therefore assume that $\pi_2 X$ is torsion. Our proof then
provides a transverse element $x \in H^*(X; \F_2)$ in even degree
whose image in $H^*(K(\pi_n X, n))$ is a transverse element of the
form $Sq^I_t u_n$ for some admissible sequence $(i_1, \dots,
i_m)$. In particular all powers $x^{2^k}$ are non-zero, which
proves (1) - (3). Moreover the elements $x, Sq^{2i_1}x, Sq^{4i_1,
2i_1} x, \dots$ are non-zero, indecomposable, and algebraically
independent because so are the corresponding images in
$H^*(K(\pi_n X, n))$. This proves (4) and~(5).
\end{proof}

\section{Cohomological finiteness conditions}
\label{section fgoverA}
The strategy we followed to analyze the integral homology of
Postnikov pieces can be applied in a more general context. We work
in this last section with simply connected $H$-spaces $X$ such
that $H^*(X; \F_2)$ is finitely generated as an algebra over the
Steenrod algebra. This section relies on the Sullivan conjecture.
As it may be considered thus as less elementary than the part
about Postnikov pieces, we have decided to postpone it till the
end of the article.

{}From the assumption on the mod $2$ cohomology, we infer by
\cite[Lemma~7.1]{deconstructing} that there exists an integer $n$
such that the module $QH^*(X; \F_2)$ of indecomposable elements
lies in the $(n-1)$-st stage of the Krull filtration for unstable
modules, \cite{MR95d:55017}. Therefore there exists by
\cite[Theorem~7.3]{deconstructing} a simply connected $H$-space $Y
= P_{B\Z/2} X$ with finite mod $2$ cohomology and a series of
principal $H$-fibrations
$$
X = X_n \xrightarrow{p_n} X_{n-1} \rightarrow \dots \rightarrow
X_1 \xrightarrow{p_1} X_0 = Y
$$
of simply connected spaces such that the homotopy fiber of $p_i$
is an Eilenberg-Mac\, Lane space $K(P_i, i)$, where $P_i$ splits
as a product of a finite direct sum $P'_i$ of cyclic groups
$\Z/2^r$ and a finite direct sum $P''_i$ of Pr\"ufer groups
$\Z_{2^\infty}$. Let us recall here that $X_k$ is obtained as the
$\Sigma^k B\Z/p$-nullification of $X$ (the above tower is
Bousfield's nullification tower, \cite{B2}). Since
$K(\Z_{2^\infty}, i)$ and $K(\Z, i+1)$ are mod $2$ equivalent, we
alter slightly the way in which the $P_i$'s are added to $Y$ in
order to work in a more familiar setting. Then we can recover $X$
from the tower
$$
X = Y_{n} \xrightarrow{q_n} Y_{n-1} \rightarrow \dots \rightarrow
Y_1 \xrightarrow{q_1} Y_0 = Y
$$
of simply connected spaces and principal $H$-fibrations, where the
homotopy fiber of $q_i$ is the product of Eilenberg-Mac\, Lane
spaces $K(P'_{i+1}, i+1) \times K(P''_i, i)$. Notice that $Q_1 =
P'_1$ is trivial because we assume that $X$ is simply connected
($Y^\wedge_2$ is therefore $2$-connected, \cite{MR23:A2201}). We
have a splitting result, just like in Proposition~\ref{proposition
factorisation du k-invariant}.

\begin{lemma}
\label{lemma factorisation du k-invariant2}
Let $X$ be a simply connected $H$-space such that $H^*(X; \F_2)$
is finitely generated as an algebra over the Steenrod algebra.
Assume that $\pi_2 X^\wedge_2$ is torsion. Then there is an
$H$-fibration
$$\xymatrix{
X \ar[r] &Y_{n-2}\ar[r] &K(P''_n \oplus P'_n,n+1) \times
K(P''_n,n). }$$
\end{lemma}

\begin{proof}
The proof is based on the Zabrodsky lemma, as in
Proposition~\ref{proposition factorisation du k-invariant}.
\end{proof}

Our next result is the analog in the present setting of
Proposition~\ref{proposition QX[n]}. Recall from
Section~\ref{section gaps} that the set $A_n$ consists in those
integers for which the $2$-adic length is strictly larger
than~$n$.

\begin{proposition}
\label{proposition QX}
Let $X$ be a simply connected $H$-space such that $H^*(X; \F_2)$
is finitely generated as an algebra over the Steenrod algebra.
There exists then integers $m$ and $N$ such that $Q^a H^*(X;\F_2)=
0 = P^a H^*(X;\F_2)$ for all $a\in A_n$ with $a \geq N$.
\end{proposition}

\begin{proof}
The integer $m$ is determined by the stage of the Krull filtration
in which $QH^*(X; \F_2)$ lives, i.e. by the degrees in which the
homotopy groups of the homotopy fiber of $X \rightarrow P_{B\Z/2}
X = Y$ are non-trivial. With the above notation, $m=n$ if $P''_n$
is trivial, and $m=n+1$ if $P''_n$ is not. The proof goes then by
induction on $m$. When $m=0$, choose $N$ to be larger than the
cohomological dimension of~$Y$. The proof of
Proposition~\ref{proposition QX[n]} goes through.
\end{proof}

\begin{lemma}
\label{lemma lowdimsplitting}
Let $X$ be a simply connected $H$-space which fits, up to
$2$-completion, in an $H$-fibration of the form
$$\xymatrix{
K(\oplus_t \Z, 2) \ar[r] & X \ar[r] & Y}
$$
where $H^*(Y; \F_2)$ is finite. Then $X$ has no homology exponent
unless the fibration splits up to $2$-completion, i.e $X
\simeq^\wedge_2 Y \times K(\oplus_t \Z, 2)$.
\end{lemma}

\begin{proof}
Let us omit the $2$-completions in the proof and write the details
of the proof when $t=1$. By the result of Hubbuck and Kane,
\cite{MR53:1582}, $\pi_3 Y$ is isomorphic to a direct sum of say
$s$ copies of~$\Z$. The map classifying the fibration factors
through $Y[3] \simeq K(\oplus_s \Z, 3)\rightarrow K(\Z, 3)$. The
$E_2$-term of the Serre spectral sequence has the form $\Z[u]
\otimes H^*(Y; \Z)$, where $u$ has degree~$2$ and the cohomology
of $Y$ is of finite dimension~$N$, and of exponent $2^a$ for some
integer~$a$. The differential $d_3(u) = x$ for some non-zero
element $x \in H^3(Y; \Z) \cong \oplus_s \Z$. Therefore $d_3(u^n)
= n x \otimes u^{n-1}$. At worst $d_3(x \otimes u^{n-1})$ is
non-zero and then hits a torsion element, of order at most $2^a$.
Hence, on the third column of the $E_4$-term, we have a group
covering $\Z/2^{n-a}$ in vertical degree $2^n$. From the
finiteness of $Y$ we see that the spectral sequence collapses at
$E_{N-3}$. An iteration of the above argument shows therefore that
the third column of the $E_\infty$-term contains a group covering
$Z/2^{n-(N-5)a}$ in vertical degree~$2^n$, for any $n \geq 1$. In
particular there is arbitrarily high torsion in $H^*(X;\Z)$.
Therefore, for $X$ to have an exponent, the fibration must split.
\end{proof}

\begin{remark}
\label{to have an exponent or not to have one?}
We point out that the preceding lemma provides simple examples of
fibrations, such as $K(\Z, 2) \rightarrow S^3\langle 3 \rangle
\rightarrow S^3$, where both the fiber and the base have an
exponent, but the total space has none.
\end{remark}

\begin{theorem}
\label{t:main2}
Let $X$ be a simply connected $H$-space of finite type such that
$H^*(X; \F_2)$ is finitely generated as an algebra over the
Steenrod algebra. Assume that $X$ admits a homology exponent. Then
$X$ is, up to $2$-completion, the product of a mod $2$ finite
$H$-space $Y$ with copies of $K(\Z, 2)$ and~$K(\Z,3)$.
\end{theorem}

\begin{proof}
We follow the proof of Theorem~\ref{t:main}. Let us thus assume
that $X$ admits a homology exponent. By killing the copies of $\Z$
in $\pi_2 X$ just like in Proposition~\ref{prop reduction}, we can
assume that $\pi_2 X^\wedge_2$ is torsion. We also see by
inspection of the tower that $\pi_2 (Y_i)^\wedge_2$ is torsion for
any $i \geq 0$. Therefore the splitting in Lemma~\ref{lemma
factorisation du k-invariant2} holds and we work with a fibration
$$\xymatrix{
K(P''_n \oplus P'_n,n) \times K(P''_n,n-1) \ar[r] & X \ar[r]
&Y_{n-2}. }
$$
If $P'_n \neq 0$, it must contain a copy of $\Z/2^r$ as direct
summand. Choose a power of the corresponding element $Sq^\xi_t
v_n$, of degree larger than the integer $N$ given in
Proposition~\ref{proposition QX}. From the Serre spectral sequence
for the above fibration we see that this provides a transverse
element in~$H^*(X;\F_2)$. Therefore $P'_n = 0$ (and so $P''_n$ is
not trivial).

If $n \geq 3$ we choose a copy of $\Z_{2^\infty}$ in $P''_n$ and a
suitable power of the corresponding element $Sq^\eta v_n$ to
detect a transverse element in $H^*(X; \F_2)$. Since we assume
that $X$ has a homology exponent, we see that $n\leq 2$, i.e. $X$
is the homotopy fiber of a map $k: Y \rightarrow K(P''_1, 2)
\times K(P''_2, 3)$. To conclude the proof we must show that this
map is trivial.

The mod $2$ cohomology of the $H$-space $Y$ is finite. Rationally
it is thus a product of odd dimensional spheres and, in
particular, $\pi_4 Y^\wedge_2$ is torsion. This implies that the
projection of $k$ on the second factor $Y \rightarrow K(P''_2, 3)$
is the trivial map, up to $2$-completion. Hence the copies of
$K(\Z^\wedge_2, 3)$ split off $X^\wedge_2$. We are left with the
analysis of a fibration $X \longrightarrow Y \longrightarrow
K(P''_1, 2)$. If the map $Y \rightarrow K(P''_1, 2)$ is not
trivial, we conclude from Lemma~\ref{lemma lowdimsplitting} that
$X$ cannot have a homology exponent. Hence, the fibration must
split and this concludes the proof.
\end{proof}


\end{document}